\documentclass[12pt]{amsart}   

\usepackage{latexsym,amscd,amsmath,verbatim,amssymb,calc,enumerate,amsxtra}
\allowdisplaybreaks[1]

\theoremstyle{plain}
\newtheorem{theorem}{Theorem}[section]
\newtheorem{corollary}[theorem]{Corollary}
\newtheorem{proposition}[theorem]{Proposition}
\newtheorem{lemma}[theorem]{Lemma}
\newtheorem*{theorem*}{Theorem}

\theoremstyle{remark}
\newtheorem{remark}[theorem]{Remark}
\newtheorem{notation}[theorem]{Notation}

\theoremstyle{definition}
\newtheorem{definition}[theorem]{Definition}

\newtheorem{example}[theorem]{Example}

\newtheorem{discussion}[theorem]{Discussion}
\newtheorem{rmd}[theorem]{Reminder}
\newtheorem*{note}{Note}

\numberwithin{equation}{theorem}

\DeclareMathOperator{\Spec}{Spec} 
\DeclareMathOperator{\Supp}{Supp}

\renewcommand{\ker}{\operatorname{Ker}}
\newcommand{\coker}{\operatorname{Coker}}
\newcommand{\image}{\operatorname{Im}}

\newcommand{\Modu}{\operatorname{Mod}}

\newcommand{\ann}{\operatorname{Ann}}

\newcommand{\spec}{\operatorname{Spec}}

\newcommand{\diag}{\operatorname{{\bf diag}}}
\newcommand{\Var}{\operatorname{Var}}

\newcommand{\un}{\operatorname{(\underline{{\it x}})-unm}} 

\renewcommand{\L}{\Lambda}
\newcommand{\m}{\mathfrak m}
\newcommand{\p}{\mathfrak p}
\newcommand{\fp}{\mathfrak p}
\newcommand{\fb}{\mathfrak b}
\newcommand{\fa}{\mathfrak a}
\newcommand{\n}{\mathfrak n}
\newcommand{\mc}{\mathcal}
\newcommand{\ida}{\mathfrak a}
\newcommand{\lra}{\longrightarrow}

\newcommand{\ul}{\underline}
\newcommand{\ol}{\overline}

\newcommand{\wh}{\widehat}

\newcommand{\bp}[1]{^{[#1]}}
\newcommand{\disappear}[1]{}

\begin{document}
\title[FROBENIUS TEST EXPONENTS]{FROBENIUS TEST EXPONENTS FOR PARAMETER IDEALS IN
GENERALIZED COHEN--MACAULAY LOCAL RINGS}
\author{CRAIG HUNEKE}
\address{Department of Mathematics,
University of Kansas, Lawrence, KS 66045, USA\\{\it Fax number}:
001-785-864-5255} \email{huneke@math.ku.edu}
\author{MORDECHAI KATZMAN}
\address{Department of Pure Mathematics,
University of Sheffield, Hicks Building, Sheffield S3 7RH, United Kingdom\\
{\it Fax number}: 0044-114-222-3769}
\email{M.Katzman@sheffield.ac.uk}
\author{RODNEY Y. SHARP}
\address{Department of Pure Mathematics,
University of Sheffield, Hicks Building, Sheffield S3 7RH, United Kingdom\\
{\it Fax number}: 0044-114-222-3769}
\email{R.Y.Sharp@sheffield.ac.uk}
\author{YONGWEI YAO}
\address{Department of Mathematics, University of Michigan, Ann Arbor,
MI 48109, USA\\ {\it Fax number}: 001-734-763-0937}
\email{ywyao@umich.edu}

\thanks{Huneke was partially supported by
the US National Science Foundation (grant number DMS-0244405);
Sharp was partially supported by the Engineering and Physical
Sciences Research Council of the United Kingdom (grant number
EP/C538803/1).}


\date{\today}


\begin{abstract}
This paper studies Frobenius powers of parameter ideals
in a commutative Noetherian local ring $R$ of prime characteristic
$p$. For a given ideal $\fa$ of $R$, there is a power $Q$ of $p$,
depending on $\fa$, such that the $Q$-th Frobenius power of the
Frobenius closure of $\fa$ is equal to the $Q$-th Frobenius power
of $\fa$. The paper addresses the question as to whether there
exists a {\em uniform\/} $Q_0$ which `works' in this context for
all parameter ideals of $R$ simultaneously.

In a recent paper, Katzman and Sharp proved that there does exists
such a uniform $Q_0$ when $R$ is Cohen--Macaulay. The purpose of
this paper is to show that such a uniform $Q_0$ exists when $R$ is
a generalized Cohen--Macaulay local ring. A variety of concepts
and techniques from commutative algebra are used, including
unconditioned strong $d$-sequences, cohomological annihilators,
modules of generalized fractions, and the
Hartshorne--Speiser--Lyubeznik Theorem employed by Katzman and
Sharp in the Cohen--Macaulay case.
\end{abstract}

\maketitle

\date{\today}

\setcounter{section}{-1}
\section{\sc Introduction}
\label{in}

This paper studies a certain type of uniform behaviour
of parameter ideals in a commutative Noetherian ring $R$ of prime
characteristic $p$.

One motivation for our work comes from the theory of test
exponents for tight closure introduced by M. Hochster and C.
Huneke in \cite[Definition 2.2]{HocHun00}. For an ideal $\fa$ of
$R$ and a non-negative integer $n$, the {\em $p^n$-th Frobenius
power\/} $\fa^{[p^n]}$ of $\fa$ is the ideal of $R$ generated by
all $p^n$-th powers of elements of $\fa$. Suppose, temporarily,
that $R$ is reduced. Recall that a {\em test element\/} for $R$ is
an element $c$ of $R$ outside all the minimal prime ideals of $R$
such that, for each ideal $\fa$ of $R$, and for $r \in R$, it is
the case that $r \in \fa^*$, the tight closure of $\fa$, if and
only if $cr^{p^n} \in \fa^{[p^n]}$ for all $n \ge 0$. It is a
result of Hochster and Huneke \cite[Theorem (6.1)(b)]{HocHun94}
that such a test element exists if $R$ is a (reduced) algebra of
finite type over an excellent local ring of characteristic $p$.

Let $c$ be a test element for $R$, and let $\fa$ be an ideal of
$R$. A {\em test exponent\/} for $c$, $\fa$ is a power $q =
p^{e_0}$ (where $e_0$ is a non-negative integer) such that if, for
an $r \in R$, we have $cr^{p^e} \in \fa^{[p^e]}$ for {\em one
single\/} $e \geq e_0$, then $r \in \fa^*$ (so that $cr^{p^n} \in
\fa^{[p^n]}$ for all $n \ge 0$). In \cite{HocHun00}, it is shown
that this concept has strong connections with the major open
problem about whether tight closure commutes with localization;
indeed, to quote Hochster and Huneke, `roughly speaking, test
exponents exist if and only if tight closure commutes with
localization'.

In a recent paper \cite{tctecpi}, R. Y. Sharp has shown that, for
a test element $c$ in a reduced equidimensional excellent local
ring $R$ (of characteristic $p$), there exists a non-negative
integer $e_0$ such that $p^{e_0}$ is a test exponent for $c$,
$\fa$ for {\em every\/} parameter ideal $\fa$ of $R$. (In such an
$R$, a parameter ideal is simply an ideal that can be generated by
part of a system of parameters.) We can think of $p^{e_0}$ as a
{\em uniform parameter test exponent\/} for $R$.

It is natural to ask whether there is an analogous result for
Frobenius closures. Return to the general situation where we
assume only that $R$ is a commutative Noetherian ring of prime
characteristic $p$. The {\em Frobenius closure $\fa^F$} of an
ideal $\fa$ of $R$ is defined by
$$ \fa ^F := \big\{ r
\in R : \mbox{there exists~} 0 \le n \in \mathbb Z \mbox{~such
that~} r^{p^n} \in \fa^{[p^n]}\big\}.
$$
This is an ideal of $R$, and so is finitely generated; therefore
there exists a power $Q_0$ of $p$ such that $(\fa^F)^{[Q_0]} =
\fa^{[Q_0]}$, and we define $Q(\fa)$ to be the smallest power of
$p$ with this property. Note that, for $r \in R$, it is the case
that $r \in \fa^F$ if and only if $r^{Q(\fa)} \in \fa^{[Q(\fa)]}$.

 In \cite[\S 0]{KS}, M. Katzman and Sharp raised
the following question: is the set $\{Q(\fb) : \fb \mbox{~is a
proper ideal of~} R\}$ of powers of $p$ bounded? In the case where
$R$ is Artinian, that is, has dimension $0$, it is easy to see
that this question has an affirmative answer, because in that case
$R/\sqrt{0}$ is a direct product of fields, and one can deduce
easily from that that, if $Q_1$ is a power of $p$ such that
$\big(\sqrt{0}\big)^{[Q_1]} = 0$, then $Q(\fb) \le Q_1$ for every
ideal $\fb$ of $R$. For this reason, we shall assume that $\dim R
> 0$ for the remainder of the paper.

H. Brenner \cite{HB} has recently shown that the answer to the
question of Katzman and Sharp, as stated above, is negative.
Nevertheless, it might not be too unreasonable to hope that, in
the case where $R$ is local, the set
$$
\{Q(\fb) : \fb \mbox{~is a parameter ideal of~} R\}
$$
is bounded. In \cite[Theorem 2.5]{KS}, Katzman and Sharp showed
that this is the case when $R$ is a Cohen--Macaulay local ring:
they showed that, then, there exists an invariant $\eta(R)$ of $R$
such that $(\fb^F)^{[p^{\eta(R)}]} = \fb^{[p^{\eta(R)}]}$ for all
parameter ideals $\fb$ of $R$. The purpose of this paper is to
prove the corresponding result when $R$ is a generalized
Cohen--Macaulay local ring, that is, when all the local cohomology
modules $H^i_{\m}(R)~(i = 0, \ldots, t-1)$ (where $\m$ denotes the
maximal ideal of $R$ and $t := \dim R$) have finite length.
Specifically, in Theorem \ref{dplus-Q} we prove the following.

\vspace{0.1in}

\noindent{\bf Theorem.} {\it Let $R$ be a generalized
Cohen--Macaulay local ring of prime characteristic $p$. Then there
exists a power $Q$ of $p$ such that $((\fb)^F)^{[Q]} = \fb^{[Q]}$
for every ideal $\fb$ of $R$ that can be generated by part of a
system of parameters of $R$.}

\vspace{0.1in}

In the Cohen--Macaulay case, the invariant $\eta (R)$ was defined
by means of the Hartshorne--Speiser--Lyubeznik Theorem (see
\cite[Theorem 1.4]{KS}) about a certain type of uniform behaviour
of a left module over the Frobenius skew polynomial ring
(associated to $R$) that is Artinian as an $R$-module: Katzman and
Sharp applied this Hartshorne--Speiser--Lyubeznik Theorem to the
top local cohomology module of a Cohen--Macaulay local ring.

In this paper, we make similar use of the
Hartshorne--Speiser--Lyubeznik Theorem, although we apply it to
all the local cohomology modules of the generalized
Cohen--Macaulay local ring $R$. We also use a variety of other
concepts and techniques from commutative algebra, including
unconditioned strong $d$-sequences and work of S. Goto and K.
Yamagishi \cite{GY} about them, filter-regular sequences,
cohomological annihilators, and modules of generalized fractions.

Other motivation for this work is provided in \cite[\S 0]{KS}.
There is no doubt in our minds that uniform behaviour of Frobenius
closures of the type established in this paper is both desirable
for its own sake and also relevant to the vigorous and ongoing
development of tight closure theory.

\section{\sc Notation and terminology}
\label{nt}

Throughout this paper, $R$ will denote a Noetherian commutative
ring with $\dim R = t > 0$, and $\fa$ will denote an ideal of $R$.
We shall use $\Var (\fa)$ to denote the {\em variety of $\fa$\/};
thus $ \Var (\fa) := \{ \fp \in \Spec (R) : \fp \supseteq \fa\}.$
We shall use $\min(R)$ to denote the set of minimal prime ideals
of $R$ and $R^{\circ}$ to denote $R \setminus \bigcup_{\fp \in
\min(R)}\fp$. The annihilator of an $R$-module $M$ will be denoted
by $\ann_R(M)$. We shall sometimes use the notation $(R,\m)$ to
indicate that $R$ is local with maximal ideal $\m$; then,
$(\widehat{R},\widehat{\m})$ will denote the $\m$-adic completion
of $R$. Also in the local case, we say that $\ul x = x_1, \dotsc,
x_l$ is a {\em system of parameters of
  $(R,\m)$} if $\sum_{i=1}^lx_iR$ is $\m$-primary and $l = t$; we say $\ul
  x$ is a {\em subsystem of parameters of $(R,\m)$} if it is a subsequence
  of a system of parameters of $R$.

\begin{notation} \label{notation}
Throughout the paper, $\ul x = x_1, x_2, \dotsc, x_l$ will denote
a sequence of $l$ elements of $R$.
\begin{enumerate}
\item We use $\mathbb N$ to denote the set of all non-negative
  integers, and $\mathbb N_+$ to denote the set of all
  positive integers.
\item The main results concern the case where $R$ has prime
characteristic $p$, but this hypothesis will only be in force when
explicitly stated; then, $q$, $q'$, $Q$, $\widetilde{Q}$ and
$Q_i~(i \in \mathbb N)$ will always denote powers of $p$ with
non-negative integer exponents. \item For integers $i \le j$, we
denote the subset $\{i,\dotsc,
  j\}$ of $\mathbb Z$ by $[i,j]$, and we agree that $[i,j] =
  \emptyset$ if $i > j$.
\item We adopt the normal convention that $a^0 = 1$ for all $a \in
R$. \item For each $\emptyset \neq \Lambda \subseteq [1,l]$, we
set
  $x_{\Lambda} := \prod_{i \in \Lambda} x_i \in R$. In case $\L = \emptyset$, we agree that
  $x_{\emptyset}^n = 1$ and $\sum_{i \in \emptyset}x_i^nR = (0)$ for all $n
  \in \mathbb N$.
\item For any $\Lambda \subseteq [1,l]$ and $n_1, \ldots, n_l \in
\mathbb N$, the sequence
$$
{\textstyle \left( \left( \left(\sum_{i \in \L} x_i^{n_i +
j}R\right) : x_{\L}^j\right) \right)_{j = 0,1,2,\ldots}}
$$
forms an ascending chain of ideals, and we denote its ultimate constant value
by $\left(\sum_{i \in \L} x_i^{n_i}R\right)^{\lim}$. Thus
$$
{\textstyle \left(\sum_{i \in \L} x_i^{n_i}R\right)^{\lim} =
\bigcup_{j\in \mathbb N} \left( \left(\sum_{i \in \L} x_i^{n_i +
j}R\right) : x_{\L}^j\right).}
$$
  In particular
  $\left(\sum_{i \in \emptyset}x_iR\right)^{\lim} = (0)$.
\item \label{xyz}
  For any $\Lambda \subseteq [1,l]$ and $n_1, \ldots, n_l \in \mathbb N$, we
  set $$
  {\textstyle \big(\sum_{i \in \L} x_i^{n_i}R\big)^{\un}
  := \big(\big(\sum_{i \in \L} x_i^{n_i}R\big) : \sum_{i \in [1,l] \setminus \L}
  x_iR\big)}$$
and refer to this as the {\em unmixed part of $\sum_{i \in \L}
x_i^{n_i}R$ relative to the sequence $\ul x = x_1, x_2, \dotsc,
x_l$}.
\end{enumerate}
\end{notation}

Next, we recall some definitions of various concepts in
commutative algebra.
\begin{definition}\label{definition} Recall that $\ida$ denotes an
ideal of $R$.

\begin{enumerate}
\item We say that $\ul x$ is an {\em $\ida$-filter regular sequence\/} if there
  exists an integer $n \in \mathbb N$ such that $\ida^n \subseteq
  \ann_R \big(\big(\big(\sum_{i=1}^{j-1}x_iR\big): x_{j}\big)/
  \sum_{i=1}^{j-1}x_iR\big)$ for all
  $j = 1,2,\dotsc,l$. It is easy to see that $\ul x$ is an $\ida$-filter regular sequence
  if and only if, for all $\fp \in \Spec (R) \setminus \Var
  (\fa)$, the natural images of $x_1, \ldots, x_l$ form a possibly
  improper regular sequence in $R_{\fp}$.
\item We say that $\ul x$ is a {\em $d$-sequence\/} if
$${\textstyle \big(\big(\sum_{i=1}^jx_iR\big) : x_{j+1}
    x_k\big) = \big(\big(\sum_{i=1}^jx_iR\big) :
    x_k\big)}$$ for all $j,k$ such that
    $0 \le j < k \le l$.
\item    If $x_1^{n_1},\dotsc, x_l^{n_l}$ form a
    $d$-sequence in any order and for any positive integers $n_1, n_2,
    \dotsc, n_l$, then we say that $\ul x$ is an
    {\em unconditioned strong $d$-sequence\/}, or a {\em $d^+$-sequence}.
\item Generalized Cohen--Macaulay local rings were studied by P.
Schenzel in \cite{Sc75} (where they were called
`quasi-Cohen--Macaulay local rings' (\cite[Definition 2]{Sc75}))
and by Schenzel, N. V. Trung and N. T. Cu\`ong in \cite{STC}.
When $(R,\m)$ is local, we say that $R$ is a {\em generalized
Cohen--Macaulay local ring\/} if
  $\bigcap_{i=0}^{t -1}\ann_R(H_{\m}^i(R))$ (recall that $t$ denotes $\dim R$)
  contains an $\m$-primary
  ideal; since all the local cohomology modules of $R$ with respect
  to $\m$ are Artinian,
  this is the case if and only if $H^i_{\m}(R)$ has finite length for all
  $i = 0, \ldots, \dim R - 1$.  Note that $R$ is a
  generalized Cohen--Macaulay local ring if and only if its completion $\widehat{R}$
  is.
\item Again when $(R,\m)$ is local, we say that $R$ is {\em is
Cohen--Macaulay on the punctured spectrum\/}
  if $R_{\p}$ is Cohen--Macaulay for every
  $\p \in \spec(R) \setminus \{\m\}$.
\end{enumerate}
\end{definition}

\begin{notation}
\label{nt.10} Suppose that $R$ has prime characteristic $p$. In
these circumstances, we shall always denote by $f:R\lra R$ the
Frobenius homomorphism, for which $f(r) = r^p$ for all $r \in R$.
We shall use the
 skew polynomial ring $R[T,f]$ associated to $R$ and $f$
in the indeterminate $T$ over $R$. Recall that $R[T,f]$ is, as a
left $R$-module, freely generated by $(T^i)_{i \in \mathbb N}$,
 and so consists
 of all polynomials $\sum_{i = 0}^n r_i T^i$, where  $n \in \mathbb N$
 and  $r_0,\ldots,r_n \in R$; however, its multiplication is subject to the
 rule
 $$
  Tr = f(r)T = r^pT \quad \mbox{~for all~} r \in R\/.
 $$
We refer to $R[T,f]$ as the {\em Frobenius skew polynomial ring
associated to $R$}.

 If $G$ is a left $R[T,f]$-module, then the set
$$\Gamma_T(G) := \left\{ g \in G : T^jg = 0 \mbox{~for some~} j
\in \mathbb N_+ \right\}$$ is an $R[T,f]$-submodule of $G$, called
the\/ {\em $T$-torsion submodule} of $G$; we say that $G$ is {\em
$T$-torsion\/} precisely when $G = \Gamma_T(G)$.

Note that $R$ itself has a natural structure as a left
$R[T,f]$-module under which $Tr = f(r)$ for all $r \in R$ (that
is, in which the action of the indeterminate $T$ on an element of
$R$ is just the same as the action of the Frobenius homomorphism).
We shall use $\!\!\phantom{i}_{R[T,f]}\!\Modu$ to denote the
category of all left $R[T,f]$-modules and $R[T,f]$-homomorphisms
between them. The category of all $R$-modules will be denoted by
$\Modu_R$.
\end{notation}

\section{Modules of generalized fractions}         
\label{mgf}

The concept of module of generalized fractions (due to
  Sharp and H. Zakeri \cite{SZ82a}) will be used in this
  paper. The construction and basic properties of these modules can be found
  in \cite{SZ82a}, but, at the request of the referee, we include in this
  section explanation of some of the main ideas.

\begin{rmd}[R. Y. Sharp and H. Zakeri
{\cite[\S 2]{SZ82a}}]\label{mgf.1} Let $k \in \mathbb N_+$. Let
$U$ be a {\em triangular subset\/} of $R^k$ \cite[2.1]{SZ82a},
that is, a non-empty subset of $R^k$ such that
\begin{enumerate}
\item whenever $(u_1, \ldots, u_k) \in U$ and $n_1, \ldots, n_k
\in \mathbb N_+$, then $(u_1^{n_1}, \ldots, u_k^{n_k}) \in U$
also; and \item whenever $(u_1, \ldots, u_k), (v_1, \ldots, v_k)
\in U$, then there exists $(w_1, \ldots, w_k) \in U$ such that $
w_i \in {\textstyle \left(\sum_{j=1}^i u_jR\right) \cap
\left(\sum_{j=1}^i v_jR\right)}$ for all $i = 1, \ldots k,$ so
that there exist $k \times k$ lower triangular matrices $\mathbf
H$ and $\mathbf K$ with entries in $R$ such that
$$
\mathbf H [u_1,\dotsc, u_k]^T = [w_1,\dotsc, w_k]^T = \mathbf K
[v_1,\dotsc, v_k]^T.
$$
(Here, $^T$ denotes matrix transpose, and $[z_1,\ldots, z_k]^T$
(for $z_1, \ldots, z_k \in R$) is to be interpreted as a $k \times
1$ column matrix in the obvious way.)
\end{enumerate}
It will be convenient for us to use $D_k(R)$ to denote the set of
$k \times k$ lower triangular matrices with entries in $R$; we use
$\det(\mathbf H)$ to denote the determinant of an $\mathbf H \in
D_k(R)$.

Let $M$ be an $R$-module. Define a relation $\sim$ on $M\times U$
as follows: for $m,n \in M$ and $(u_1, \ldots, u_k), (v_1, \ldots,
v_k) \in U$, write $(m, (u_1, \ldots, u_k)) \sim (n, (v_1, \ldots,
v_k))$ precisely when there exist $(w_1, \ldots,w_k) \in U$ and
$\mathbf H, \mathbf K \in D_k(R)$ such that
$$
\mathbf H [u_1,\ldots, u_k]^T = [w_1,\ldots, w_k]^T = \mathbf K
[v_1,\ldots, v_k]^T
$$
and $\det(\mathbf H)m - \det(\mathbf K)n \in {\textstyle
\sum_{j=1}^{k-1} w_jM.} $

Then $\sim$ is an equivalence relation; for $m \in M$ and $(u_1,
\ldots, u_k) \in U$, we denote the equivalence class of $(m, (u_1,
\ldots, u_k))$ by the `generalized fraction'
$$ \frac{m}{(u_1,\ldots, u_k)}.$$
The set of all equivalence classes of $\sim$ is an $R$-module,
called the {\em module of generalized fractions of $M$ with
respect to $U$}, under operations for which, for $m,n \in M$ and
$(u_1, \ldots, u_k), (v_1, \ldots, v_k) \in U$,
$$
\frac{m}{(u_1,\ldots, u_k)} + \frac{n}{(v_1,\ldots, v_k)} =
\frac{\det(\mathbf H)m + \det(\mathbf K)n}{(w_1,\ldots, w_k)}
$$
for {\em any\/} choice of $(w_1, \ldots, w_k) \in U$ and $\mathbf
H, \mathbf K \in D_k(R)$ such that
$$
\mathbf H [u_1,\ldots, u_k]^T = [w_1,\ldots, w_k]^T = \mathbf K
[v_1,\ldots, v_k]^T,
$$
and, for $r \in R$,
$$
r\frac{m}{(u_1,\ldots, u_k)} = \frac{rm}{(u_1,\ldots, u_k)}.
$$
This module of generalized fractions is denoted by $U^{-k}M$.
\end{rmd}

\begin{remark} \label{mgf.1a}
Use the notation of \ref{mgf.1}. It is worth bearing in mind, when
one is calculating with generalized fractions in $U^{-k}M$, that,
whenever $(u_1, \ldots, u_k), (y_1, \ldots, y_k) \in U$ and
$\mathbf L \in D_k(R)$ are such that $\mathbf L[u_1,\ldots, u_k]^T
= [y_1,\ldots, y_k]^T$, then, for all $m \in M$,
$$
\frac{m}{(u_1,\ldots, u_k)} = \frac{\det(\mathbf L)m}{(y_1,\ldots,
y_k)} \quad \mbox{in $U^{-k}M$}.
$$
This sometimes permits us to change a denominator to one that is,
in some sense, more convenient.
\end{remark}

\begin{example} \label{mgf.2}
For our
  given sequence $\ul x = x_1, x_2, \dotsc, x_l$ and any $i \in \mathbb N_+$,
  we set
  $$U(\ul x)_i := \{(x_1^{n_1}, \dotsc, x_i^{n_i}) : \mbox{for some~} j \in
  [0,i], n_1,
  \dotsc,n_j \in \mathbb N_+, n_{j+1}= \dotsb = n_i = 0\},$$
  where $x_k$ is interpreted as $1$ when $k \ge l+1$.  Then $U(\ul x)_i$ is a
  triangular subset of $R^i$
  for each $i \in \mathbb N_+$.
\end{example}

\begin{discussion} \label{mgf.3} We use the notation of \ref{mgf.1}.

\begin{enumerate}

\item The construction of a module of generalized fractions can be
viewed as a generalization of ordinary fraction formation in
commutative algebra: see \cite[3.1]{SZ82a}.

\item For $r_1, \ldots, r_k \in R$, we shall denote the diagonal
$k \times k$ matrix in $D_k(R)$ whose diagonal entries are $r_1,
\ldots, r_k$ by $\diag(r_1, \ldots, r_k)$. Notice that
$$
\det (\diag(r_1, \ldots, r_k)) = r_1\cdots r_k,
$$
and that, for non-negative integers $n_1, \ldots, n_k, m_1,
\ldots,m_k$, we have
$$
\diag(r_1^{m_1}, \ldots, r_k^{m_k})[ r_1^{n_1}, \ldots,
r_k^{n_k}]^T = [r_1^{m_1+n_1}, \ldots, r_k^{m_k+n_k}]^T.
$$

\item The comments in (ii) can be useful.  For example, if $(u_1,
\ldots, u_k), (w_1, \ldots, w_k)$ $\in U$ are such that there
exist $\mathbf H, \mathbf K \in D_k(R)$ with
$$
\mathbf H [u_1, \ldots, u_k]^T = [w_1, \ldots, w_k]^T = \mathbf K
[u_1, \ldots, u_k]^T,
$$
then
$$
\mathbf D\mathbf H [u_1, \ldots, u_k]^T = [w_1^2, \ldots, w_k^2]^T
= \mathbf D\mathbf K [u_1, \ldots, u_k]^T,
$$
where $\mathbf D := \diag (w_1, \ldots, w_k)$, and it turns out
(see \cite[Lemma 2.3]{SZ82a}) that $\det(\mathbf D\mathbf H) -
\det(\mathbf D\mathbf K) \in \sum_{i=1}^{k-1} w_i^2R$. This can be
helpful when one has rather little information about $\mathbf H$
but full knowledge of $\mathbf K$.

To illustrate this, suppose that it is known, for some $r \in R$,
that $\det(\mathbf H)r \in \sum_{i=1}^{k-1}w_iR$; then $
\det(\mathbf D\mathbf H)r = w_1\cdots w_{k}\det(\mathbf H)r \in
\sum_{i=1}^{k-1} w_i^2R, $ and the above considerations show that
$$\det(\mathbf D\mathbf K)r
= \det(\mathbf D\mathbf H)r - \left(\det(\mathbf D\mathbf H) -
\det(\mathbf D\mathbf K)\right)r \in {\textstyle \sum_{i=1}^{k-1}
w_i^2R}.
$$

\item The equivalence relation $\sim$ of \ref{mgf.1} is such that,
whenever $(u_1, \ldots, u_k) \in U$ and the element $m$ of $M$
actually belongs to $\sum_{i=1}^{k-1}u_iM$, then
$$
\frac{m}{(u_1, \ldots, u_k)} = 0 \quad \mbox{in $U^{-k}M$}.
$$

\item Note that, by (iv), if $k \ge 2$ and $(v_1, \ldots,
v_{k-2},1,1) \in U$, then
$$
\frac{n}{(v_1, \ldots,v_{k-2},1,1)} = 0 \quad \mbox{in $U^{-k}M$,
for all $n \in M$}.
$$

\item In the case where $U$ consists entirely of possibly improper
regular sequences on $M$, there is the following converse of (iv),
proved by Sharp and Zakeri in \cite[Theorem 3.15]{SZ82c}: if
$(u_1, \ldots, u_k) \in U$ and $m \in M$ are such that
$$
\frac{m}{(u_1, \ldots, u_k)} = 0,
$$
then $m \in \sum_{i=1}^{k-1}u_iM$.
\end{enumerate}
\end{discussion}

\begin{discussion}
\label{mgf.4} A {\em chain of triangular subsets on $R$\/} is a
family $\mc U
  := (U_i)_{i=1}^{\infty}$ such that
  \begin{enumerate}
  \item $U_i$ is a triangular subset of $R^i$ for all $i \in\mathbb N_+$;
  \item whenever $(u_1, \ldots,u_{i}) \in U_i$ (for any $i \in \mathbb N_+$),
  then $(u_1, \ldots,u_{i},1) \in U_{i+1}$;
  \item whenever $(u_1, \ldots,u_{i}) \in U_i$ with $i > 1$, then $(u_1,
  \ldots,u_{i-1}) \in U_{i-1}$;
  and
  \item $(1) \in U_1$.
  \end{enumerate}
Given such a family $\mc U$ and an $R$-module $M$, we can
construct a complex $\mc C(\mc U, M)$ of modules of generalized
fractions
\[
0 \xrightarrow{e^{-1}} M \xrightarrow{e^{0}} U_1^{-1}M
\xrightarrow{e^{1}} \dotsb \xrightarrow{e^{i-1}} U_i^{-i}M
\xrightarrow{e^{i}}  U_{i+1}^{-i-1}M \to \dotsb,
\]
in which $e^0(m) = \frac {m}{(1)} \in U_1^{-1}M$ for all $m \in M$
and, for $i \in \mathbb N_+$,
$$e^i \left(\frac {m}{(u_1, \dotsc, u_i)} \right)
= \frac {m}{(u_1, \dotsc, u_i,1)} \in U_{i+1}^{-i-1}M$$ for all $m
\in M$ and $(u_1, \dotsc, u_i) \in
 U_i.$ The comment in \ref{mgf.3}(v) shows that
 $\mc C(\mc U, M)$ is indeed a complex.

  The Exactness Theorem of Sharp--Zakeri (\cite[Theorem
 3.3]{SZ82c}, but see L. O'Carroll \cite[Theorem 3.1]{O'Car83} for
 a subsequent shorter proof) states that the complex $\mc C(\mc U,
 M)$ is exact if and only if, for all $i \in \mathbb N_+$, each element of
 $U_i$ is a possibly improper regular sequence on $M$.
\end{discussion}

\begin{example} \label{mgf.5}
For our
  sequence $\ul x = x_1, x_2, \dotsc, x_l$, the family $\mc U(\ul x)
  := (U(\ul x)_i)_{i=1}^{\infty}$, where the $U(\ul x)_i~(i\in \mathbb N_+)$ are as defined in
  \ref{mgf.2}, is a chain of triangular subsets on
  $R$. In particular, we can construct the complex of modules of
  generalized fractions $\mc C(\mc U(\ul x),R)$, which we shall write as
\[
0 \xrightarrow{e^{-1}} R \xrightarrow{e^{0}} U(\ul x)_1^{-1}R
\xrightarrow{e^{1}} \dotsb \xrightarrow{e^{i-1}} U(\ul x)_i^{-i}R
\xrightarrow{e^{i}}  U(\ul x)_{i+1}^{-i-1}R \to \dotsb.
\]
  Note that $U(\ul x)_i^{-i}R =
0$ whenever $i \ge l+2$. When working with a cohomology module
$H^i(\mc C(\mc U(\ul x),R)) = \ker e^i/\image e^{i-1}$ (where $i
\in \mathbb N$) of the complex $\mc C(\mc U(\ul x),R)$, we shall
use `$\left[\phantom{x_1,x_2}\right]$' to denote natural images,
in this cohomology module, of elements of $\ker e^i$.

  In the special case in which $(R,\m)$ is local, $l = \dim R = t$
  and $x_1, \ldots, x_l$ is a system of parameters of $R$, it
  follows from the Exactness Theorem mentioned in \ref{mgf.4} that
  $R$ is Cohen--Macaulay if and only if $\mc C(\mc U(\ul x),R)$ is
  exact. In \cite[Theorems 2.4 and 2.5]{SZ85}, Sharp and Zakeri proved that
  $(R,\m)$ is a generalized Cohen--Macaulay local ring if and
  only if there is a power of $\m$ that annihilates all the
  cohomology modules of the complex $\mc C(\mc U(\ul x),R)$, and
  that when this is the case, the $i$th cohomology module $H^i(\mc C(\mc U(\ul x),R))$ of the
  complex $\mc C(\mc U(\ul x),R)$ is isomorphic to $H^i_{\m}(R)$
  for all $i = 0, \ldots, \dim R - 1 = l-1$. The latter result is
  relevant to this paper, although in \S \ref{pr} we provide
  a refinement for use in the case where $R$ has prime
  characteristic $p$.

\end{example}

Some of the applications of modules of generalized fractions can
be found in
  \cite{SZ82b}, \cite{SZ82c}, \cite{SZ85}, \cite{O'Car83}, \cite{40} and
  \cite{KSZ}. This paper provides some more.

\section{\sc Background results}
\label{br}

In this section, we shall collect together some observations and
results, with appropriate references, that concern topics
mentioned in \S \ref{nt} and which we plan to use.

\begin{remark} \label{dplus}
Suppose that $\ul x = x_1, \dotsc, x_l$ is an unconditioned strong
$d$-sequence in $R$.
\begin{enumerate}
\item It is immediate from the definitions that, for any $n_1,
\dotsc, n_l \in \mathbb N_+$, $\L \subsetneq [1,l]$
  and $j \in [1,l] \setminus \L$, we have $${\textstyle \left(\left(\sum_{i \in \L} x_i^{n_i}R
  \right)   : x_j^{n_j}\right) = \left(\left(\sum_{i \in \L} x_i^{n_i}R\right) : x_j\right)
  = \left(\sum_{i \in \L} x_i^{n_i}R\right)^{\un}.}$$
\item Therefore, if $\ida$ is an ideal of $R$ such that $\ida
\subseteq \sqrt {\sum_{i=1}^l x_iR}$, then every permutation of
$x_1, \dotsc, x_l$ is an $\ida$-filter regular
  sequence, because $\sum_{i=1}^l x_iR$ annihilates
  ${\textstyle \big(\sum_{i \in \L} x_iR : \sum_{i \in [1,l] \setminus\L} x_iR\big)
  \big/\left(\sum_{i \in \L} x_iR\right)}$ for all $\L \subsetneq [1,l]$.
\end{enumerate}
\end{remark}

\begin{theorem} [P. Schenzel
{\cite[Satz 2.4.2]{Sch}}] \label{br.1} Suppose that $(R,\m)$ is
local; recall that $\dim(R) = t$. Then, for all systems of
parameters $a_1, \ldots, a_t$ of $R$, the ideal $\prod_{i=0}^{t-1}
\ann_R(H_{\m}^i(R))$ of $R$ annihilates all the $R$-modules
$$
{\textstyle \big(\big(\sum_{i=1}^ja_iR\big):
a_{j+1}\big)/\big(\sum_{i=1}^ja_iR\big)  \quad  (j = 0,\ldots,
t-1)}.
$$
\end{theorem}

In \cite{GO}, S. Goto and T. Ogawa proved that, in a generalized
Cohen--Macaulay local ring $(R,\m)$, there exists a positive
integer $h$ such that every system of parameters of $R$ contained
in $\m^h$ is a $d$-sequence. The following corollary (of
Schenzel's Theorem \ref{br.1}) is a variation on that theme.

\begin{corollary} \label{br.1a} Suppose that $(R,\m)$ is local
with $\dim(R) = t$, and that $\ul x = x_1, \dotsc, x_l$ is a
subsystem of parameters of $R$ such that $\sum_{i=1}^lx_iR
\subseteq \prod_{i=0}^{t-1} \ann_R(H_{\m}^i(R))$. Then $\ul x$ is
an unconditioned strong $d$-sequence.
\end{corollary}

\begin{proof} Let $n_1, \ldots, n_l$ be positive integers, and let
$j\in \mathbb N$ be such that $0 \le j < l$. By Schenzel's Theorem
\ref{br.1}, the $R$-module $\big(\big(\sum_{i=1}^jx_i^{n_i}R\big):
x_{j+1}^{n_{j+1}}\big)/\big(\sum_{i=1}^jx_i^{n_i}R\big)$ is
annihilated by $\sum_{i=1}^lx_iR$ and so, in particular, by
$x_{j+1}$. Hence
$$
{\textstyle \big(\big(\sum_{i=1}^jx_i^{n_i}R\big):
x_{j+1}^{n_{j+1}}\big) = \big(\big(\sum_{i=1}^jx_i^{n_i}R\big):
x_{j+1}\big).}$$ Moreover, the hypotheses on $x_1, \ldots, x_l$ do
not depend on the order in which this sequence is written.

Therefore, if $j,k$ are such that $0 \le j < k \le l$, then the
fact that $x_k$ annihilates
$$
{\textstyle \big(\big(\sum_{i=1}^jx_i^{n_i}R\big):
x_{j+1}^{n_{j+1}}\big)/\big(\sum_{i=1}^jx_i^{n_i}R\big)}
$$
ensures that $ {\textstyle \big(\big(\sum_{i=1}^jx_i^{n_i}R\big):
x_{j+1}^{n_{j+1}}x_k^{n_k}\big) \subseteq
\big(\big(\sum_{i=1}^jx_i^{n_i}R\big): x_{k}^{n_{k}+1}\big),} $
and the preceding paragraph shows that $$ {\textstyle
\big(\big(\sum_{i=1}^jx_i^{n_i}R\big): x_{k}^{n_{k}+1}\big) =
\big(\big(\sum_{i=1}^jx_i^{n_i}R\big): x_{k}^{n_k}\big).} $$
Therefore $x_1^{n_1}, \ldots, x_l^{n_l}$ is a $d$-sequence. Since
the hypotheses on $x_1, \ldots, x_l$ do not depend on the order in
which this sequence is written, we see that $\ul x$ is an
unconditioned strong $d$-sequence.
\end{proof}

The next proposition is an extension of a well-known result. For
an explanation of what it means to say that a local ring is
formally catenary, and for a proof of L. J. Ratliff's Theorem that
a universally catenary local ring is formally catenary, the reader
is referred to \cite[p.\ 252]{HM}.

\begin{proposition} \label{br.2}
Assume that $(R,\m)$ is a formally catenary local ring all of
whose formal fibres are Cohen--Macaulay. (These hypotheses
 would be satisfied if $R$ was an excellent local ring.)

Suppose that $R$
 is equidimensional and Cohen--Macaulay on the punctured spectrum.  Then
  $R$ is a generalized Cohen--Macaulay local ring.
\end{proposition}

\begin{proof} Since $R$ is equidimensional and formally catenary,
$\widehat{R}$ is equidimensional; the
hypothesis concerning the formal fibres ensures that $\widehat{R}$
is Cohen--Macaulay on the punctured spectrum. Thus one can assume
that $R$ is complete, and in that case the claim follows from work
of P. Schenzel, N. V. Trung and N. T. Cu\`ong \cite[(2.5) and
(3.8)]{STC}.
\end{proof}

The next two results are entirely due to S. Goto and K. Yamagishi
\cite{GY}. Unfortunately, as far as we are aware, \cite{GY} only
exists as a preprint that has been circulating informally for more
than 15 years, without formal publication. T. Kawasaki included a
proof of \cite[Lemma~2.2]{GY} in \cite[Theorem A.1]{Ka}; as we
have not been able to find a formally published proof of
\cite[Theorem~2.3]{GY}, we have included one below.

\begin{theorem} [S. Goto and K. Yamagishi
{\cite[Lemma~2.2]{GY}}]
  \label{colon} {\rm (See also \cite[Theorem A.1]{Ka}.)}
Suppose that $\ul x = x_1, \dotsc,x_{l}$ is an unconditioned
strong $d$-sequence in $R$. Then, for each $\Delta \subsetneq
[1,l]$ and each $j \in [1,l]\setminus \Delta$, and for all
positive integers $n_1,\dotsc, n_{l}$, we have \begin{align*}
{\textstyle \big(\sum_{i \in \Delta}x_i^{n_i}R\big)^{\un}} & =
{\textstyle \big(\big(\sum_{i \in \Delta}x_i^{n_i}R\big) :
x_j\big)}\\& = {\textstyle \sum_{\Lambda \subseteq \Delta}
\big(\prod_{i \in \Lambda} x_i^{n_i - 1}\big) \big(\sum_{i \in
\L}x_iR\big)^{\un}.}\end{align*}
\end{theorem}

\begin{proof} This result, originally due to Goto and Yamagishi, is proved in
\cite[Theorem~A.1]{Ka} in the case where $n_1,\dotsc, n_{l} \ge
2$, and one can check that that proof works for all choices of
positive integers $n_1,\dotsc, n_{l}$.
\end{proof}

\begin{theorem} [S. Goto and K. Yamagishi
{\cite[Theorem~2.3]{GY}}] \label{GY} Suppose that $\ul x = x_1,
\dotsc, x_l$ is an unconditioned strong $d$-sequence in $R$ and
let $n_1, n_2,\dotsc, n_l \in \mathbb N_+$ be any positive
integers.
\begin{enumerate}
\item When $l=1$,
  $${\textstyle (x_1^{n_1}R)^{\lim} = x_1^{n_1}R + \bigcup_{j \in \mathbb N}(0:
  x_1^{j}) = x_1^{n_1}R + (0 : x_1).}$$
\item When $l \ge 2$,
\begin{align*}
{\textstyle \big(\sum_{i=1}^lx_i^{n_i}R\big)^{\lim}} & =
{\textstyle \sum_{i=1}^l \big(\big(\sum_{j \in [1,l]\setminus
\{i\}} x_j^{n_j} R \big):_R x_i \big)}\\ & = {\textstyle
\sum_{i=1}^l \big(\sum_{j \in [1,l]\setminus \{i\}} x_j^{n_j}
R\big)^{\un}.}
\end{align*}
\end{enumerate}
\end{theorem}

\begin{proof} (Goto--Yamagishi.) We use induction on $l$; the
result is easy when $l = 1$, and so we suppose that $l \ge 2$ and
that the result has been proved for smaller values of $l$.

Let $a \in \big(\sum_{i=1}^lx_i^{n_i}R\big)^{\lim}$, so that,
without loss of generality, there exists $m \in \mathbb N_+$ such
that $x_{[1,l]}^ma \in \sum_{i=1}^lx_i^{n_i+m}R$. Thus there exist
$b \in \sum_{i=1}^{l-1}x_i^{n_i+m}R$ and $c \in R$ such that $
x_{[1,l]}^ma = b + x_l^{n_l+m}c$, and so
$$
{\textstyle x_{[1,l-1]}^ma - x_l^{n_l}c \in
\big(\big(\sum_{i=1}^{l-1}x_i^{n_i+m}R\big):x_l^m\big) =
\big(\big(\sum_{i=1}^{l-1}x_i^{n_i+m}R\big):x_l\big).}
$$
Therefore, by \ref{colon}, we have \[ {\textstyle x_{[1,l-1]}^ma -
x_l^{n_l}c = \sum_{\L \subseteq [1,l-1]}\prod_{i \in\L}x_i^{n_i
+m-1}h_{\L}}, \tag{$\ddagger$} \] where $h_{\L} \in
\big(\sum_{i\in\L}x_iR\big)^{\un}$ for all $\L \subseteq [1,l-1]$.
However, for each $\Lambda \subsetneq [1,l-1]$ and $j \in [1,l-1]
\setminus \L$, we have
$${\textstyle \big(\big(\sum_{i\in\L}x_iR\big):x_l\big) =
\big(\sum_{i\in\L}x_iR\big)^{\un} =
\big(\big(\sum_{i\in\L}x_iR\big):x_j\big)}$$ by \ref{dplus}(i), so
that $x_{[1,l-1]}h_{\L} \in \sum_{i\in\L}x_i^2R$. We now multiply
both sides of equation $(\ddagger)$ by $x_{[1,l-1]}$ to obtain
that
$$
{\textstyle x_{[1,l-1]}^{m+1}\left(a -
\prod_{i=1}^{l-1}x_i^{n_i-1}h_{[1,l-1]}\right) \in
\sum_{i=1}^{l-1}x_i^{n_i+m+1}R + x_l^{n_l}R.}
$$
Since the natural images of $x_1, \ldots, x_{l-1}$ in the ring
$R/x_l^{n_l}R$ form an unconditioned strong $d$-sequence in that
ring, it follows from the inductive hypothesis that
\begin{align*}
{\textstyle a - \prod_{i=1}^{l-1}x_i^{n_i-1}h_{[1,l-1]}} & \in
{\textstyle \sum_{i=1}^{l-1} \big(\big(\sum_{j \in
[1,l-1]\setminus \{i\}} x_j^{n_j} R + x_l^{n_l}R \big):_R x_i
\big) + x_l^{n_l}R + x_1^{n_1}R}
\\ & = {\textstyle \sum_{i=1}^{l-1} \big(\big(\sum_{j \in [1,l-1]\setminus
\{i\}} x_j^{n_j} R + x_l^{n_l}R \big):_R x_i \big) + x_1^{n_1}R.}
\end{align*}
(The presence of the ideal $x_1^{n_1}R$ on the right hand side
ensures that the argument applies to the case where $l = 2$.)
Since $\big(\sum_{i=1}^{l-1}x_iR\big)^{\un} =
\big(\big(\sum_{i=1}^{l-1}x_iR\big):x_l\big),$ we see that
$${\textstyle \prod_{i=1}^{l-1}x_i^{n_i-1}h_{[1,l-1]} \in
\big(\big(\sum_{j=1}^{l-1} x_j^{n_j} R \big):_R x_l \big)}$$ (even
in the case where $l = 2$) and so $a \in \sum_{i=1}^{l}
\big(\big(\sum_{j \in [1,l]\setminus \{i\}} x_j^{n_j} R \big):_R
x_i \big).$ Therefore $${\textstyle
\big(\sum_{i=1}^{l}x_i^{n_i}R\big)^{\lim} \subseteq \sum_{i=1}^l
\big(\big(\sum_{j \in [1,l]\setminus \{i\}} x_j^{n_j} R \big):_R
x_i \big) = \sum_{i=1}^l \big(\sum_{j \in [1,l]\setminus \{i\}}
x_j^{n_j} R\big)^{\un},}$$ and the reverse inclusion is easy.
\end{proof}

The following corollary is immediate from \ref{colon} and
\ref{GY}.

\begin{corollary}
\label{br.21} Suppose that $\ul x = x_1, \dotsc, x_l$ is an
unconditioned strong $d$-sequence in $R$ (with $l \ge 1$) and let
$n_1, n_2,\dotsc, n_l \in \mathbb N_+$ be any positive integers.
Then
\begin{enumerate}
\item $x_{[1,l]}\big(\sum_{i=1}^{l}x_i^{n_i}R\big)^{\lim}
\subseteq \sum_{i=1}^{l}x_i^{n_i+1}R$; and \item
$\big(\sum_{i=1}^{l}x_i^{n_i}R\big)^{\lim} = \sum_{\Lambda
\subsetneq [1,l]} \big(\prod_{i \in \Lambda} x_i^{n_i - 1}\big)
\left(\sum_{i \in \Lambda} x_iR\right)^{\un}$ when $l \ge 2$.
\end{enumerate}
\end{corollary}

\begin{theorem} \label{Hu}  Let $\ul x = x_1,x_2,\dotsc,x_l$ be a
$d$-sequence in $R$.

\begin{enumerate}
\item {\rm (C. Huneke \cite[Proposition~2.1]{Hu}.)} Then
$$ {\textstyle \big(\big(\sum_{i=1}^rx_iR\big): x_{r+1}\big)
\bigcap \big(\sum_{i=1}^{l}x_iR\big) = \sum_{i=1}^rx_iR} \quad
\text{~for all~} r = 0, \ldots, l-1.
$$
 \item If $\ul x = x_1,x_2,\dotsc,x_l$ is an
  unconditioned strong $d$-sequence in $R$, then
\[
{\textstyle \big(\sum_{i \in \L} x_i^{n_i}R \big)^{\un} \bigcap
\big(\sum_{i=1}^l x_i^{n_i} R\big) = \sum_{i \in \L} x_i^{n_i}R}
\]
for all positive integers $n_1, n_2,\dotsc, n_l$ and each $\L
\subsetneq [1,l]$.
\end{enumerate}
\end{theorem}

\begin{proof} (ii) This is immediate from Huneke's result quoted in part
(i), because $${\textstyle \big(\sum_{i \in \L} x_i^{n_i}R
\big)^{\un} = \big( \big(\sum_{i \in \L} x_i^{n_i}R \big) : x_j
\big) = \big( \big(\sum_{i \in \L} x_i^{n_i}R \big) : x_j^{n_j}
\big)} \quad \text{ for~} j \in [1,l] \setminus \L\text{:}$$ see
\ref{dplus}(i).
\end{proof}

As in \cite{KS}, use will be made of the following extension, due
to G. Lyubeznik, of a result of R. Hartshorne and R. Speiser. It
shows that, when $R$ is local and of prime characteristic $p$, a
$T$-torsion left $R[T,f]$-module which is Artinian (that is,
`cofinite' in the terminology of Hartshorne and Speiser) as an
$R$-module exhibits a certain uniformity of behaviour.

\begin{theorem} [G. Lyubeznik {\cite[Proposition 4.4]{Ly}}]
\label{HSL}  {\rm (Compare Hartshorne--Speiser \cite[Proposition
1.11]{HS}.)} Suppose that $(R,\m)$ is local and of prime
characteristic $p$, and let $G$ be a left $R[T,f]$-module which is
Artinian as an $R$-module. Then there exists $e \in \mathbb N$
such that $T^e\Gamma_T(G) = 0$.
\end{theorem}

Hartshorne and Speiser first proved this result in the case where
$R$ is local and contains its residue field which is perfect.
Lyubeznik applied his theory of $F$-modules to obtain the result
without restriction on the local ring $R$ of characteristic $p$. A
short proof of the theorem, in the generality achieved by
Lyubeznik, is provided in \cite{HSLonly}.

\begin{lemma} [Katzman--Sharp {\cite[Lemma 3.5]{KS}}]
\label{br.3} Suppose that $R$ has prime characteristic $p$.
\begin{enumerate}
\item Let $n \in \mathbb N_+$ and let $U$ be a triangular subset
of $R^n$. Then the module of generalized fractions $U^{-n}R$ has a
structure as left $R[T,f]$-module with
$$
T\left(\frac{r}{(u_1, \ldots,u_n)}\right) = \frac{r^p}{(u_1^p,
\ldots,u_n^p)} \quad \mbox{~for all~} r \in R \mbox{~and~} (u_1,
\ldots,u_n) \in U.
$$
\item It follows easily that the complex $\mc C(\mc U(\ul x), R)$
of modules of generalized fractions of\/ {\rm \ref{mgf.5}} is a
complex of left $R[T,f]$-modules and $R[T,f]$-homomorphisms; hence
all its cohomology modules $H^i(\mc C(\mc U(\ul x),R))~(i \in
\mathbb N)$ have natural structures as left $R[T,f]$-modules.
\end{enumerate}
\end{lemma}

\section{Preparatory results}
\label{pr}

Most of the results in this section concern the case where $R$ has
prime characteristic $p$, but the first two do not.

\begin{lemma}\label{gen-frac} Suppose that $(R,\m)$ is local, and consider the complex
of modules of generalized fractions $\mc C(\mc U(\ul x), R)$ of\/
{\rm \ref{mgf.5}.} Let $r$ be an integer such that $0 \le r < l$.
(In the case where $r = 0$, a generalized fraction such as
$h/(x_1, \dotsc, x_r)$ (where $h \in R$) is to be interpreted
simply as $h$.)
\begin{enumerate}
\item If $h \in \big(\big(\sum_{i=1}^rx_iR\big): x_{r+1}\big)$,
then
$$\frac{h}{(x_1,
  \dotsc, x_r)} \in \ker e^r, \quad \mbox{so that~} \left[\frac{h}{(x_1,
  \dotsc, x_r)}\right] \in H^r(\mc C(\mc U(\ul x), R)).$$ (The
  notation `$\left[\phantom{x_1x_2}\right]$' is explained in\/
  {\rm \ref{mgf.5}}.)
\item Let $n \in \mathbb N_+$ and $h \in R$. Then
    $$
  \frac{h}{(x_1^n, \dotsc,
  x_r^n)} \in \image e^{r-1}$$
  if and only if $h \in (\sum_{i=1}^rx_i^nR\big)^{\lim}$.
\end{enumerate}
\end{lemma}

\begin{proof}
When $r=0$, all the claims are easy.  We therefore omit the proofs
in that case and assume that $1 \le r < l$.

(i) Since $x_{r+1}h \in \sum_{i=1}^r x_iR$, it is immediate from
\ref{mgf.1a} and \ref{mgf.3}(iv) that
$$e^r \left(\frac{h}{(x_1, \dotsc,
x_r)}\right) = \frac{h}{(x_1, \dotsc, x_r,1)} = \frac{x_{r+1}
h}{(x_1, \dotsc, x_r,x_{r+1})} =0 \in U(\ul x)_{r+1}^{-r-1}R.$$

(ii) ($\Leftarrow$) Assume that $h \in
(\sum_{i=1}^rx_i^nR\big)^{\lim}$. Thus there exists $m \in \mathbb
N$ such that $x_{[1,r]}^m h \in \sum_{i=1}^rx_i^{n+m}R$. Thus we
can write $x_{[1,r]}^m h = \sum_{i=1}^r s_{i} x_i^{n+m}$ for some
$s_1, \ldots, s_r \in R$. Then, in $U(\ul x)_r^{-r}R$, we have
\begin{align*}
\frac{h}{(x_1^n, \dotsc,x_{r-1}^n, x_r^n)} & = \frac{x_{[1,r]}^m
h}{(x_1^{n+m}, \dotsc, x_{r-1}^{n+m}, x_r^{n+m})}
 = \frac{\sum_{i = 1}^r s_{i} x_i^{n+m}}
{(x_1^{n+m}, \dotsc, x_{r-1}^{n+m}, x_r^{n+m})}\\
 & = \frac{s_{r} x_r^{n+m}}{(x_1^{n+m}, \dotsc, x_{r-1}^{n+m},
 x_r^{n+m})}\quad \mbox{~(on use of \ref{mgf.3}(iv))~}\\
 &= \frac{s_{r}}{(x_1^{n+m}, \dotsc, x_{r-1}^{n+m},
 1)}\\ & = \begin{cases}
 e^{r-1}(s_r)                          &    \text{if } r = 1, \\
 e^{r-1}\left(\frac{s_{r}} {(x_1^{n+m}, \dotsc, x_{r-1}^{n+m}
 )}\right)             &    \text{if } r \ge 2.
\end{cases}
\end{align*}

($\Rightarrow$) By \ref{mgf.1a}, we can write, in $U(\ul
x)_r^{-r}R$,
$$
\frac{h}{(x_1^n, \dotsc, x_r^n)} = e^{r-1}
\left(\frac{g}{(x_1^{m+n}, x_2^{m+n}, \dotsc, x_{r-1}^{m+n})}
\right) = \frac{x_r^{m+n} g}{(x_1^{m+n}, x_2^{m+n}, \dotsc,
x_r^{m+n})} $$ for some $m \in \mathbb N$ and $g \in R$. Thus
$$\frac{x_{[1,r]}^{m} h
  - x_r^{m+n}g}{(x_1^{m+n}, x_2^{m+n}, \dotsc, x_r^{m+n})} = 0
  \quad \mbox{~in~} U(\ul x)_r^{-r}
R.$$ By the definition of modules of generalized fractions (see
\ref{mgf.1}), there exist $u \in \mathbb N$ and $\mathbf H \in
D_r(R)$ such that
\[
\mathbf H [x_1^{m+n},\dotsc, x_r^{m+n}]^T = [x_1^{m+n+u},\dotsc,
x_r^{m+n+u}]^T\] and $\det (\mathbf H) (x_{[1,r]}^{m} h  -
x_r^{m+n}g) \in \sum_{i=1}^{r-1}x_i^{m+n+u}R$. Since
$$\diag(x_1^u, \ldots,x_r^u)[x_1^{m+n},\dotsc, x_r^{m+n}]^T =
[x_1^{m+n+u},\dotsc, x_r^{m+n+u}]^T,$$ we can use the method of
\ref{mgf.3}(iii) to see that
$$
x_{[1,r]}^{m + n + u}x_{[1,r]}^u (x_{[1,r]}^{m} h  - x_r^{m+n}g)
\in {\textstyle \sum_{i=1}^{r-1}x_i^{2m+2n+2u}R}.
$$
Consequently, $x_{[1,r]}^{2m + n + 2u} h \in
\sum_{i=1}^{r-1}x_i^{2m+2n+2u}R + x_r^{2m+2n+2u}R$; this implies
that
\[
{\textstyle h \in \big(\big(\sum_{i=1}^{r}x_i^{2m+2n+2u}R\big) :
x_{[1,r]}^{2m + n + 2u}\big) \subseteq
\big(\sum_{i=1}^rx_i^{n}R\big)^{\lim}},
\]
as required.
\end{proof}

It is an immediate consequence of Corollary \ref{br.1a} that a
generalized Cohen--Macaulay local ring has a system of parameters
that is an unconditioned strong $d$-sequence. We now use modules
of generalized fractions to establish a converse of this.

\begin{theorem} \label{br.1b} Suppose that $(R,\m)$ is local;
recall that $\dim(R) = t$. Then the following statements are
equivalent:
\begin{enumerate}
\item $R$ is generalized Cohen--Macaulay; \item there exists $h
\in \mathbb N_+$ such that $y_1^k, \ldots, y_t^k$ is an
undconditioned strong $d$-sequence, for every system of parameters
$y_1, \ldots, y_t$ of $R$ and every $k \ge h$; \item there exists
a system of parameters of $R$ which is an unconditioned strong
$d$-sequence.
\end{enumerate}
\end{theorem}

\begin{proof} (i) $\Rightarrow$ (ii) This is immediate from Corollary
\ref{br.1a} and the definition of generalized Cohen--Macaulay
local ring: just choose $h$ so that $\m^h \subseteq
\bigcap_{i=0}^{t -1}\ann_R(H_{\m}^i(R))$.

(ii) $\Rightarrow$ (iii) This is clear.

(iii) $\Rightarrow$ (i) Take $l = t$ and $x_1, \ldots, x_t$ to be
a system of parameters of $R$ that is an unconditioned strong
$d$-sequence. By \ref{dplus}(ii), every permutation of $x_1,
\ldots, x_t$ is an $\m$-filter regular sequence. In order to show
that $R$ is generalized Cohen--Macaulay, it is enough, by
symmetry, to show that $x_{i+1}H^i_{\m}(R) = 0$ for all $i =
0,\ldots, t-1$.

We now apply \cite[Corollary 2.3 and Theorem 2.4]{SZ85} to the
complex $\mc C(\mc U(\ul x),R)$ of \ref{mgf.5}. Note that, for all
$i \in \mathbb N_+$, every element of $U(\ul x)_i$ is an
$\m$-filter regular sequence. The cited results from \cite{SZ85}
therefore show that
$$
H^i(\mc C(\mc U(\ul x),R)) = \ker e^i/\image e^{i-1} \cong
H^i_{\m}(R) \quad \mbox{for all $i = 0, \ldots,t-1$}.
$$
It is therefore enough for us to show that, for an $i \in
\{0,\ldots,t-1\}$, we have $x_{i+1}\ker e^i \subseteq \image
e^{i-1}$. We shall deal here with the case where $i > 0$, and
leave to the reader the (easy) modification for the case where $i
= 0$.

Let $\alpha \in \ker e^i$. By \ref{mgf.1a}, we can write
$$ \alpha = \frac{r}{(x_1^n,\ldots,x_i^n)} \quad \mbox{for some $r
\in R$ and $n \in \mathbb N_+$}. $$ Therefore
$$
\frac{x_{i+1}^nr}{(x_1^n,\ldots,x_i^n,x_{i+1}^n)} =
\frac{r}{(x_1^n,\ldots,x_i^n,1)} = e^i(\alpha) = 0.
$$
By \ref{mgf.1} and \ref{mgf.1a}, this means that there exist $v
\in \mathbb N_+$ and $\mathbf H \in D_{i+1}(R)$ such that
$$
\mathbf H [x_1^n, \ldots, x_{i+1}^n]^T = [x_1^{n+v}, \ldots,
x_{i+1}^{n+v}]^T \quad \mbox{and} \quad \det(\mathbf H)x_{i+1}^nr
\in {\textstyle \sum_{j=1}^ix_j^{n+v}R}.
$$
Since $\diag(x_1^v, \ldots, x_{i+1}^v)[x_1^n, \ldots, x_{i+1}^n]^T
= [x_1^{n+v}, \ldots, x_{i+1}^{n+v}]^T$, we can use the technique
of \ref{mgf.3}(iii) to see that $x_{[1,i+1]}^{n+2v}x_{i+1}^nr \in
\sum_{j=1}^ix_j^{2n+2v}R$. Therefore, since $x_1, \ldots,x_t$ is
an unconditioned strong $d$-sequence,
$$
x_{[1,i]}^{n+2v}r  \in \left({\textstyle \sum_{j=1}^ix_j^{2n+2v}R}
: x_{i+1}^{2n+2v}\right) = \left({\textstyle
\sum_{j=1}^ix_j^{2n+2v}R} : x_{i+1}\right). $$ Hence $x_{i+1}r \in
\left({\textstyle \sum_{j=1}^ix_j^{2n+2v}R} :
x_{[1,i]}^{n+2v}\right) \subseteq \left({\textstyle
\sum_{j=1}^ix_j^{n}R}\right)^{\lim}$, so that
$$
x_{i+1}\alpha = x_{i+1}\frac{r}{(x_1^n,\ldots,x_i^n)} =
\frac{x_{i+1}r}{(x_1^n,\ldots,x_i^n)} \in \image e^{i-1}
$$
by Lemma \ref{gen-frac}(ii).
\end{proof}

It is well known that, when $R$ has prime characteristic $p$, each
local cohomology module $H^i_{\fa}(R)$, where $i \in \mathbb N$,
has a natural structure as a left $R[T,f]$-module. A detailed
explanation is given in \cite[2.1]{KS}, and the argument there can
easily be modified to show that, if $M$ is an arbitrary left
$R[T,f]$-module, then $H^i_{\fa}(M)$ (formed by regarding $M$ as
an $R$-module by restriction of scalars) inherits a natural
structure as a left $R[T,f]$-module. However, in this paper, we
are going to use the following rather stronger statement.

\begin{proposition}
\label{pr.11} Suppose that $R$ has prime characteristic $p$. Then
$\big( H^i_{\mathfrak{a}} \big)_{i \in \mathbb N}$ is a negative
strongly connected sequence of covariant functors from
$\!\!\phantom{i}_{R[T,f]}\!\Modu$ to itself.
\end{proposition}

\begin{note} We identify $H^0_{\fa}$ with the $\fa$-torsion functor $\Gamma_{\fa}$
in the natural way. If $M$ is a left $R[T,f]$-module, then
$\Gamma_{\fa}(M)$ is an $R[T,f]$-submodule of $M$. It should be
noted from the proof below that this $R[T,f]$-module structure on
$\Gamma_{\fa}(M)$ is exactly the same as the natural left
$R[T,f]$-module structure on $H^0_{\fa}(M) = \Gamma_{\fa}(M)$
provided by the proposition.
\end{note}

\begin{proof} As this proof relies on the Independence
Theorem for local cohomology (see \cite[4.2.1]{BS}), we shall use
notation similar to that employed in \cite[\S 4.2]{BS}. Let
$\lceil : \Modu_R \longrightarrow \Modu_R$ denote the functor
obtained from restriction of scalars using the Frobenius
homomorphism $f$: thus, if $Y$ is an $R$-module, then $Y\lceil$
denotes $Y$ considered as an $R$-module via $f$.

Let $M$ be a left $R[T,f]$-module. The map $\tau_M : M
\longrightarrow M\lceil$ defined by $\tau_M(m) = Tm$ is an
$R$-module homomorphism. Consequently, for each $i \in \mathbb N$,
there is an induced $R$-homomorphism $H^i_{\fa}(\tau_M) :
H^i_{\fa}(M) \longrightarrow H^i_{\fa}(M\lceil)$.

Let
\[
\Theta = (\theta^i)_{i \in \mathbb N} : \left(H^i
_{\mathfrak{a}}(\: {\scriptscriptstyle \bullet} \: \lceil)
\right)_{i \in \mathbb N} \stackrel{\cong}{\longrightarrow}
\left(H^i _{\mathfrak{a}^{[p]}}(\: {\scriptscriptstyle \bullet}
\:) \lceil \right)_{i \in \mathbb N}
\]
be the isomorphism of negative (strongly) connected sequences of
covariant functors (from $\Modu_R$ to $\Modu_R$) that is inverse
to the one given in the Independence Theorem for local cohomology,
in the form in which it is stated in \cite[4.2.1]{BS}. Thus
$\theta^0$ is the identity natural equivalence. Since $\fa$ and
$\fa^{[p]}$ have the same radical, $H^i_{\fa}$ and
$H^i_{\fa^{[p]}}$ are the same functor, for each $i \in \mathbb
N$.

Consider the $\mathbb Z$-endomorphism $\xi^i_M := \theta^i_{M}
\circ H^i_{\fa}(\tau_M) : H^i_{\fa}(M) \lra H^i_{\fa}(M)\lceil$.
We now modify the argument of  \cite[2.1]{KS} and use \cite[Lemma
1.3]{KS} to show that $H^i_{\fa}(M)$ has a natural structure as a
left $R[T,f]$-module in which $T\gamma = \xi^i_M(\gamma)$ for all
$\gamma \in H^i_{\fa}(M)$.

Furthermore, if $\beta : M \longrightarrow N$ is an
$R[T,f]$-homomorphism of left $R[T,f]$-modules, then
\[
\begin{picture}(300,75)(-150,-25)
\put(-40,40){\makebox(0,0){$ M $}} \put(40,40){\makebox(0,0){$
M\lceil$}} \put(-25,40){\vector(1,0){50}}
\put(-36,10){\makebox(0,0)[l]{$^{ \beta }$}}
\put(0,46){\makebox(0,0){$^{ \tau_M }$}}
\put(0,-14){\makebox(0,0){$^{ \tau_N }$}}
\put(44,10){\makebox(0,0)[l]{$^{ \beta\lceil }$}}
\put(-40,-20){\makebox(0,0){$ N $}} \put(40,-20){\makebox(0,0){$
N\lceil$}} \put(-25,-20){\vector(1,0){50}}
\put(-40,30){\vector(0,-1){40}} \put(40,30){\vector(0,-1){40}}
\end{picture}
\]
is a commutative diagram of $R$-homomorphisms, so that, for $i \in
\mathbb N$, the diagram
\[
\begin{picture}(300,75)(-150,-25)
\put(-80,40){\makebox(0,0){$ H^i_{\mathfrak{a}}(M) $}}
\put(-40,46){\makebox(0,0){$^{ H^i_{\mathfrak{a}}(\tau_M) }$}}
\put(0,40){\makebox(0,0){$ H^i_{\mathfrak{a}}(M\lceil) $}}
\put(40,46){\makebox(0,0){$^{ \theta^i_{M} }$}}
\put(80,40){\makebox(0,0){$ H^i_{\mathfrak{a}}(M)\lceil $}}
\put(-58,40){\vector(1,0){36}} \put(22,40){\vector(1,0){36}}
\put(-76,10){\makebox(0,0)[l]{$^{ H^i_{\mathfrak{a}}(\beta) }$}}
\put(4,10){\makebox(0,0)[l]{$^{ H^i_{\mathfrak{a}}(\beta\lceil)
}$}} \put(84,10){\makebox(0,0)[l]{$^{
H^i_{\mathfrak{a}}(\beta)\lceil }$}} \put(-80,-20){\makebox(0,0){$
H^i_{\mathfrak{a}}(N) $}} \put(-40,-14){\makebox(0,0){$^{
H^i_{\mathfrak{a}}(\tau_N) }$}} \put(0,-20){\makebox(0,0){$
H^i_{\mathfrak{a}}(N\lceil) $}} \put(40,-14){\makebox(0,0){$^{
\theta^i_{N} }$}} \put(80,-20){\makebox(0,0){$
H^i_{\mathfrak{a}}(N)\lceil $}} \put(-58,-20){\vector(1,0){36}}
\put(22,-20){\vector(1,0){36}} \put(-80,30){\vector(0,-1){40}}
\put(0,30){\vector(0,-1){40}} \put(80,30){\vector(0,-1){40}}
\end{picture}
\]
also commutes. This means that, when $H^i_{\fa}(M)$ and
$H^i_{\fa}(N)$ are given their natural structures as left
$R[T,f]$-modules, as in the preceding paragraph, then the
$R$-homomorphism $H^i_{\fa}(\beta) : H^i_{\fa}(M) \longrightarrow
H^i_{\fa}(N)$ is an $R[T,f]$-homomorphism. In this way,
$H^i_{\fa}$ becomes a functor from
$\!\!\phantom{i}_{R[T,f]}\!\Modu$ to itself.

Next, whenever $0 \longrightarrow L
\stackrel{\alpha}{\longrightarrow} M
\stackrel{\beta}{\longrightarrow} N \longrightarrow 0$ is an exact
sequence of left $R[T,f]$-modules and $R[T,f]$-homomorphisms, the
diagram
\[
\begin{picture}(300,65)(-150,-15)
\put(-98,40){\makebox(0,0){$ 0 $}} \put(-50,40){\makebox(0,0){$ L
$}} \put(0,40){\makebox(0,0){$                      M $}}
\put(50,40){\makebox(0,0){$                     N $}}
\put(98,40){\makebox(0,0){$ 0 $}} \put(-90,40){\vector(1,0){30}}
\put(-40,40){\vector(1,0){30}} \put(10,40){\vector(1,0){30}}
\put(60,40){\vector(1,0){30}} \put(-46,15){\makebox(0,0)[1,l]{$^{
\tau_L }$}} \put(4,15){\makebox(0,0)[l]{$^{ \tau_M }$}}
\put(54,15){\makebox(0,0)[l]{$^{                \tau_N }$}}
\put(-98,-10){\makebox(0,0){$ 0 $}} \put(-50,-10){\makebox(0,0){$
L\lceil $}} \put(0,-10){\makebox(0,0){$ M\lceil $}}
\put(50,-10){\makebox(0,0){$                      N\lceil $}}
\put(98,-10){\makebox(0,0){$ 0 $}} \put(-25,46){\makebox(0,0){$^{
\alpha }$}} \put(25,46){\makebox(0,0){$^{                   \beta
}$}} \put(-25,-4){\makebox(0,0){$^{                   \alpha\lceil
}$}} \put(25,-4){\makebox(0,0){$^{                    \beta\lceil
}$}} \put(-90,-10){\vector(1,0){30}}
\put(-40,-10){\vector(1,0){30}} \put(10,-10){\vector(1,0){30}}
\put(60,-10){\vector(1,0){30}} \put(-50,30){\vector(0,-1){30}}
\put(0,30){\vector(0,-1){30}} \put(50,30){\vector(0,-1){30}}
\end{picture}
\]
of $R$-modules and $R$-homomorphisms commutes, and so the vertical
maps induce a morphism of the long exact sequence of local
cohomology modules of the upper sequence into that for the lower
sequence. It follows from this (and properties of the isomorphism
$\Theta$ of connected sequences) that the connecting
$R$-homomorphisms
$$
H^i_{\fa}(N) \longrightarrow H^{i+1}_{\fa}(L) \quad (i \in \mathbb
N)
$$ are all homomorphisms of left $R[T,f]$-modules. Hence the long
exact sequence of local cohomology $R$-modules induced by $0
\longrightarrow L \stackrel{\alpha}{\longrightarrow} M
\stackrel{\beta}{\longrightarrow} N \longrightarrow 0$ is actually
a long exact sequence of left $R[T,f]$-modules and
$R[T,f]$-homomorphisms.

Everything else needed for completion of the proof is now
straightforward.
\end{proof}

Our next result can be viewed as a strengthening, in the
particular case where $R$ has prime characteristic $p$, of special
cases of \cite[Theorem (2.4)]{SZ85} and of results of K.
Khashyarmanesh, Sh.\ Salarian and H. Zakeri in \cite[Theorem 1.2
and Consequences 1.3(i)]{KSZ} (which refer for proof to the proof
of \cite[Theorem (2.4)]{SZ85}).

\begin{theorem}\label{iso} Suppose that $(R,\m)$ is local and of
prime characteristic $p$, and that $\ul x = x_1, \dotsc, x_l$ is
an\/ $\m$-filter regular sequence of elements of\/ $\m$. Consider
the complex $\mc C(\mc U(\ul x),R)$ of modules of generalized
fractions of\/ {\rm \ref{mgf.5}}, and note that, by\/ {\rm
\ref{br.3}(ii)}, this is a complex of left $R[T,f]$-modules and
$R[T,f]$-homomorphisms. Then there are isomorphisms of
$R[T,f]$-modules
$$H^i(\mc C(\mc U(\ul x),R)) \cong H_{\m}^i(R) \quad \mbox{~for
all~} i = 0, \ldots, l-1,$$ where the $H_{\m}^i(R)$ are considered
as left $R[T,f]$-modules in the natural way described in\/ {\rm
\ref{pr.11}}.
\end{theorem}

\begin{proof} First, it follows from \cite[3.2]{SZ82a} and
\cite[2.2]{SZ82b} that $H^j_{\m}\left(U({\ul x})_i^{-i}R\right) =
0$ for all $i = 1, \ldots, l$ and all $j \ge 0$. Second, one can
use \ref{definition}(i) and the Exactness Theorem for complexes of
modules of generalized fractions (see \ref{mgf.4}) (in conjunction
with \cite[Proposition 2.1]{40}) to see that $\Supp \left( H^i(\mc
C(\mc U(\ul x),R))\right) \subseteq \{\m\}$ for all $i \geq 0$.

With these observations, the theorem can be proved by an obvious
modification of the argument used to prove \cite[Theorem
(2.4)]{SZ85}, provided one notes that all the sequences $$0
\longrightarrow \ker e^0 \longrightarrow R \longrightarrow \image
e^0 \longrightarrow 0,$$
$$
0 \longrightarrow \ker e^i \longrightarrow U(\ul x)_i^{-i}R
\longrightarrow \image e^i \longrightarrow 0 \quad (1 \le i \le
l),
$$
$$
0 \longrightarrow \image e^{i-1} \longrightarrow U(\ul x)_i^{-i}R
\longrightarrow \coker e^{i-1} \longrightarrow 0 \quad (1 \le i
\le l)
$$
and
$$
0 \longrightarrow \image e^{i-1} \longrightarrow \ker e^i
\longrightarrow \ker e^i/\image e^{i-1} \longrightarrow 0 \quad (1
\le i \le l)
$$
are exact sequences of left $R[T,f]$-modules and
$R[T,f]$-homomorphisms, so that, by Proposition \ref{pr.11}, all
the isomorphisms of local cohomology modules that they induce are
$R[T,f]$-isomorphisms.
\end{proof}

\begin{corollary} \label{fixedQ}
Suppose that $(R,\m)$ is local and of prime characteristic $p$;
recall that $\dim(R) = t$. Since the local cohomology modules
$H^i_{\m}(R)$ are left $R[T,f]$-modules that are Artinian as
$R$-modules, it follows from the Hartshorne--Speiser--Lyubeznik
Theorem\/ {\rm \ref{HSL}} that there exists $e_1 \in \mathbb N$
such that $$T^{e_1}\Gamma_T(H^i_{\m}(R)) = 0 \quad \text{~for
all~} i = 0, \ldots, t-1.$$

Set $Q_1 := p^{e_1}$. Then $Q_1$ has the following property:
whenever $\ul x = x_1, x_2, \dotsc, x_t$ is a system of parameters
of $R$ that is also an $\m$-filter regular sequence, and whenever
$r$ is an integer with $0 \le r < t$ and $h \in
\big(\sum_{i=1}^rx_iR\big)^{\un}$ is such that $ h^q \in
\big(\sum_{i=1}^rx_i^qR\big)^{\lim}$ for some $q$, then $h^{Q_1}
\in \big(\sum_{i=1}^rx_i^{Q_1}R\big)^{\lim}.$
\end{corollary}

\begin{proof} In the case where $r = 0$, a generalized fraction such as
$h/(x_1,  \dotsc, x_r)$ (where $h \in R$) is to be interpreted
simply as $h$.  Let $\ul x = x_1, x_2, \dotsc, x_t$, $r$ and $h$
be as in the statement of the corollary.

Consider the complex of modules of generalized fractions $\mc
C(\mc U(\ul x), R)$ of \ref{mgf.5}. Since $h \in
\big(\sum_{i=1}^rx_iR\big)^{\un} \subseteq
\big(\big(\sum_{i=1}^rx_iR\big): x_{r+1}\big)$ by
\ref{notation}(vii), it follows from Lemma \ref{gen-frac}(i) that
$$\frac{h}{(x_1,
  \dotsc, x_r)} \in \ker e^r, \quad \mbox{so that~} \left[\frac{h}{(x_1,
  \dotsc, x_r)}\right] \in H^r(\mc C(\mc U(\ul x), R)).$$
Since $h^q \in \big(\sum_{i=1}^rx_i^qR\big)^{\lim}$, it follows
from Lemma \ref{gen-frac}(ii) that
$$
\left[\frac{h}{(x_1,
  \dotsc, x_r)}\right] \in \Gamma_T \left(H^r(\mc C(\mc U(\ul x),
  R))\right).
  $$
Now $H^r(\mc C(\mc U(\ul x), R)) \cong H^r_{\m}(R)$ as left
$R[T,f]$-modules, by Theorem \ref{iso}. Therefore
$$
\left[\frac{h^{Q_1}}{(x_1^{Q_1},
  \dotsc, x_r^{Q_1})}\right] = T^{e_1}\left[\frac{h}{(x_1,
  \dotsc, x_r)}\right] = 0, \quad \mbox{~so that~}
  \frac{h^{Q_1}}{(x_1^{Q_1},
  \dotsc, x_r^{Q_1})} \in \image e^{r-1}.
$$
Therefore $h^{Q_1} \in \big(\sum_{i=1}^rx_i^{Q_1}R\big)^{\lim}$ by
Lemma \ref{gen-frac}(ii).
\end{proof}

\begin{proposition} \label{lim}
Suppose that $(R,\m)$ is local and of prime characteristic $p$;
recall that $\dim(R) = t$. Then there exists $Q_2$ such that, for
each system of parameters $\ul x = x_1, \dotsc, x_t$ of $R$, we
have $\big(\big(\sum_{i=1}^{t}x_iR\big)^F\big)\bp {Q_2} \subseteq
\big(\sum_{i=1}^{t}x_i^{Q_2}R\big)^F \bigcap
\big(\sum_{i=1}^{t}x_i^{Q_2}R\big)^{\lim}$. (Note that
$\sum_{i=1}^{t}x_i^{Q_2}R =
\big(\sum_{i=1}^{t}x_iR\big)^{[Q_2]}$.)
\end{proposition}

\begin{proof} Our intention is to apply the
Hartshorne--Speiser--Lyubeznik Theorem~\ref{HSL} to the top local
cohomology module $H^t_{\m}(R)$ of $R$. Recall that $H^t_{\m}(R)$
can be realized as the $t$-th cohomology module of the \u{C}ech
complex of $R$ with respect to $x_1, \ldots, x_t$. Thus
$H^t_{\m}(R)$ can be represented as the residue class module of
$R_{x_1 \cdots x_t}$ modulo the image, under the \u{C}ech
`differentiation' map, of $\bigoplus_{i=1}^tR_{x_1 \cdots
x_{i-1}x_{i+1}\cdots x_t}$. See \cite[\S 5.1]{BS}. We use
`$\left[\phantom{=} \right]$' to denote natural images of elements
of $R_{x_1 \cdots x_t}$ in this residue class module.

Recall also (from, for example, \cite[2.3]{KS}) that the natural
left $R[T,f]$-module structure on $H^t_{\m}(R)$ is such that
$$
T\left[\frac{r}{(x_1\cdots x_t)^n}\right] =
\left[\frac{r^p}{(x_1\cdots x_t)^{np}}\right] \quad \mbox{~for
all~} r \in R \mbox{~and~} n \in \mathbb N.
$$

Since $H^t_{\m}(R)$ is an Artinian $R$-module, it follows from the
Hartshorne--Speiser--Lyubeznik Theorem \ref{HSL} that there exists
$e_2 \in \mathbb N$ such that
$T^{e_2}\Gamma_T\left(H^t_{\m}(R)\right) = 0$. Set $Q_2 =
p^{e_2}$.

Let $a \in \big(\sum_{i=1}^{t}x_iR\big)^F$, so that there exists
$Q = p^e$ such that $a^Q \in \big(\sum_{i=1}^{t}x_iR\big)^{[Q]} =
\sum_{i=1}^{t}x_i^QR$. Thus, in $H^t_{\m}(R)$, we have
$$
T^e\left[\frac{a}{x_1\cdots x_t}\right] =
\left[\frac{a^Q}{(x_1\cdots x_t)^Q}\right] = 0, \quad \text{so
that} \left[\frac{a^{Q_2}}{(x_1\cdots x_t)^{Q_2}}\right] =
T^{e_2}\left[\frac{a}{x_1\cdots x_t}\right] = 0.
$$
By \cite[(2.3)(i)]{KS}, this means that $a^{Q_2} \in
\big(\sum_{i=1}^{t}x_i^{Q_2}R\big)^{\lim}$. It is easy to check
that $\big(\big(\sum_{i=1}^{t}x_iR\big)^F\big)\bp {Q_2} \subseteq
\big(\sum_{i=1}^{t}x_i^{Q_2}R\big)^F$, and so the proof is
complete.
\end{proof}

\begin{lemma}\label{complete}
Suppose that $(R,\m)$ is local and of prime characteristic $p$. If
the ideal $\fa$ of $R$ satisfies $((\fa \wh R)^F)\bp{Q} = (\fa \wh
R)\bp{Q},$ then $(\fa^F)\bp{Q} = \fa\bp{Q}.$
\end{lemma}

\begin{proof} Let
$r \in \fa^F$; then $r \in (\fa \wh R)^F$ in the ring $\wh R$.
Thus $r^{Q} \in (\fa \wh R)\bp{Q}\cap R = (\fa\bp Q \wh R)\cap R$,
and the latter ideal is just $\fa\bp{Q}$ because $\wh R$ is a
faithfully flat extension of $R$.
\end{proof}

\begin{lemma}\label{nil}
Suppose that $R$ is of prime characteristic $p$. If the ideals
$\fa$ and $\n$ of $R$ satisfy $\n \bp {Q'} = 0$ and
$(((\ida+\n)/\n)^F)\bp{\widetilde{Q}} =
((\ida+\n)/\n)\bp{\widetilde{Q}}$ (in $R/\n$), then
$(\ida^F)\bp{Q'\widetilde{Q}} = \ida\bp{Q'\widetilde{Q}}$.
\end{lemma}

\begin{proof}
Let $r \in \ida^F$; then $r+ \n \in ((\ida+\n)/\n)^F$ in $R/\n$.
Thus $(r + \n)^{\widetilde{Q}} \in
((\ida+\n)/\n)\bp{\widetilde{Q}}$; that is, $r^{\widetilde{Q}} \in
\ida\bp{\widetilde{Q}} + \n$. Consequently, $r^{Q'\widetilde{Q}}
\in \ida\bp{Q'\widetilde{Q}}$.
\end{proof}

\section{The main results}         
\label{mr}

Throughout this section, we assume that $R$ has prime characteristic
$p$.

\begin{proposition}\label{mr.1}
Suppose that $R$ is of prime characteristic $p$. Assume that $R$
is semi-local or that the integral closure of $R/\sqrt 0$ in its
total ring of fractions is module-finite over $R/\sqrt 0$ (this
would be the case if $R$ was excellent). Then there exists $Q_3$
such that $((xR)^F) \bp {Q_3} = (xR)\bp {Q_3}$ for all $x \in
R^{\circ} := R \setminus \bigcup_{\fp
  \in \min(R)} \fp$.
\end{proposition}

\begin{proof}
In case $R$ is semi-local, as everything involved commutes with
localization at the finitely many maximal ideals of $R$, we can
assume that $(R,\m)$ is local. Then, by Lemma~\ref{complete}, we
can further assume that $(R,\m)$ is complete and hence excellent.

Thus, also by Lemma~\ref{nil}, we can assume that $R$ is reduced
and that $\ol R$ is module-finite over $R$, where $\ol R$ is the
integral closure of $R$ in its total fraction ring
$(R^{\circ})^{-1} R$. Consider $\left(\ol R \cap
R^{1/q}\right)_{q=1}^{\infty}$, which forms an ascending chain of
$R$-submodules of $\ol R$. As $\ol R$ is module-finite over $R$,
there exists $Q$ such that $\ol R \cap R^{1/q} = \ol R \cap
R^{1/Q}$ for all $q \ge Q$.

For any $x \in R^{\circ}$ and any $y \in (xR)^F$, there exists $q$
such that $y^q = a x^q$ for some $a \in R$. This means that
$(y/x)^q = a/1$ in $(R^{\circ})^{-1} R$, and this implies that
$y/x \in  \ol R \cap R^{1/q}$. By our choice of $Q$, we get $y/x
\in \ol R \cap R^{1/Q}$. Thus $(y/x)^Q = b/1$ for some $b \in R$
and hence $y^Q = b x^Q \in (xR)\bp Q$.
\end{proof}

The next theorem is the main result of this paper. Recall from
Theorem \ref{br.1b} that a local ring is generalized
Cohen--Macaulay if and only if it has a system of parameters that
is an unconditioned strong $d$-sequence.

\begin{theorem}\label{dplus-Q}
Suppose that $(R,\m)$ is a generalized Cohen--Macaulay local ring
of prime characteristic $p$; recall that $\dim(R) = t > 0$. Then
there exists $Q$ such that ${\textstyle \big(\big(\sum_{i=1}^j
x_iR\big)^F\big)\bp {Q} = \big(\sum_{i=1}^j x_iR\big)\bp {Q}}$ for
all subsystems of parameters $x_1, \dotsc, x_j$ of $R$.
\end{theorem}

\begin{proof}
In view of Proposition \ref{mr.1}, we can assume that $t \ge 2$.

In the first part of the proof, we are going to show that there
exists $Q_0$ such that $\big(\big(\sum_{i=1}^tx_iR\big)^F\big)\bp
{Q_0} = \big(\sum_{i=1}^tx_iR\big)\bp {Q_0}$ for all systems of
parameters $\ul x = x_1, \dotsc, x_t$ of $R$ that are
unconditioned strong $d$-sequences.

Let $Q_1$ be as in Corollary~\ref{fixedQ}. Also, by
Proposition~\ref{lim}, there exists $Q_2$ such that
$\big(\big(\sum_{i=1}^{t}x_iR\big)^F\big)\bp {Q_2} \subseteq
\big(\sum_{i=1}^{t}x_i^{Q_2}R\big)^F \bigcap
\big(\sum_{i=1}^{t}x_i^{Q_2}R\big)^{\lim}$ for all systems of
parameters $\ul x = x_1, \dotsc, x_t$ of $R$. Set $Q_0 = pQ_1
Q_2$. We are going to show that
$\big(\big(\sum_{i=1}^tx_iR\big)^F\big)\bp {Q_0} =
\big(\sum_{i=1}^tx_iR\big)\bp {Q_0}$ for all systems of parameters
$\ul x = x_1, \dotsc, x_t $ of $R$ that are unconditioned strong
$d$-sequences. Notice that
\begin{align*}
{\textstyle \big(\big(\sum_{i=1}^{t}x_iR\big)^F\big)\bp {pQ_1Q_2}}
& = {\textstyle \big(\big(\big(\sum_{i=1}^{t}x_iR\big)^F\big)\bp
{Q_2}\big)\bp {pQ_1}}\\ & \subseteq {\textstyle
\big(\big(\sum_{i=1}^{t}x_i^{Q_2}R\big)^F \bigcap
\big(\sum_{i=1}^{t}x_i^{Q_2}R\big)^{\lim}\big)\bp {pQ_1}}
\end{align*}
by Proposition~\ref{lim}. Therefore, it suffices to prove that
$$
{\textstyle \big(\big(\sum_{i=1}^{t}x_i^{Q_2}R\big)^F \bigcap
\big(\sum_{i=1}^{t}x_i^{Q_2}R\big)^{\lim}\big)\bp {pQ_1} \subseteq
\big(\sum_{i=1}^{t}x_i^{Q_2}R\big)\bp {pQ_1}}
$$
for each system of parameters $\ul x = x_1, \dotsc, x_t$ such that
$x_1^{Q_2}, \ldots, x_t^{Q_2}$ is an unconditioned strong
$d$-sequence, and so it is enough to prove that
$${\textstyle \big(\big(\sum_{i=1}^{t}x_iR\big)^F \bigcap
\big(\sum_{i=1}^{t}x_iR\big)^{\lim}\big)\bp {pQ_1} \subseteq
\big(\sum_{i=1}^{t}x_iR\big)\bp {pQ_1}}$$ for all systems of
parameters $\ul x = x_1, \dotsc, x_t$ of $R$ that are
unconditioned strong $d$-sequences. We therefore fix a typical
such $\ul x = x_1, \dotsc, x_t$. Notice that $\ul x$ is also an
$\m$-filter regular sequence in any order, by Remark
\ref{dplus}(ii).

Let $y \in \big(\sum_{i=1}^{t}x_iR\big)^F \cap
\big(\sum_{i=1}^{t}x_iR\big)^{\lim}$. Then there exists $q'=pq$
such that $y^{q'} \in \big(\sum_{i=1}^{t}x_iR\big)\bp {q'}$, and
without loss of generality we can assume that $q \ge
\max\{p,Q_1\}$.  We see, from Corollary~\ref{br.21}(ii), that
$y^{p} \in \big(\big(\sum_{i=1}^{t}x_iR\big)^{\lim}\big)\bp {p}
\subseteq \big(\sum_{i=1}^{t}x_i^pR\big)^{\lim} = \sum_{\L
\subsetneq [1,t]} x_{\L}^{p-1} \big(\sum_{i \in \L}
x_iR\big)^{\un}$. We can therefore write $y^p = \sum_{\L
\subsetneq [1,t]} x_{\L}^{p-1} h_{\L}$ with $h_{\L} \in
\left(\sum_{i \in \L} x_iR\right)^{\un}$ for all $\L \subsetneq
[1,t]$. Consequently, we have
\[
{\textstyle \sum_{\L \subsetneq [1,t]} x_{\L}^{pq-q} h_{\L}^q =
y^{pq} \in \big(\sum_{i=1}^tx_iR\big)\bp{pq},} \tag{$*_0$}
\]
in which $h_{\L}^q \in \big(\big(\sum_{i \in \L}
x_iR\big)^{\un}\big)\bp q \subseteq \big(\sum_{i \in \L}
x_i^qR\big)^{\un}$ for all $\L \subsetneq [1,t]$ (in view of
\ref{dplus}(i)).

The immediate goal is to show that $y^{pQ_1} \in
\big(\sum_{i=1}^tx_iR\big)\bp{pQ_1}$. To this end, as $y^{pQ_1} =
\sum_{\L \subsetneq [1,t]} x_{\L}^{pQ_1-Q_1} h_{\L}^{Q_1}$, it is
enough for us to prove that
\[
{\textstyle x_{\L}^{pQ_1-Q_1} h_{\L}^{Q_1} \in \sum_{i \in
\L}x_i^{pQ_1}R} \quad \text{~for all $\L$ such that $\L \subsetneq
[1,t]$.} \tag{$\dagger$}
\]

We now prove ($\dagger$) by induction on $|\L|$, the cardinality
of $\L$. When $|\L| = 0$, we have $\L= \emptyset$; we recall our
conventions that $\sum_{i \in \emptyset}x_i^{pQ_1}R = (0)$ and
$x_{\emptyset} = 1$. When we consider $R$ as a left
$R[T,f]$-module as in \ref{nt.10}, the $R$-submodule
$\Gamma_{\m}(R)$ is actually a $T$-torsion $R[T,f]$-submodule.
Since $\big(\sum_{i \in \emptyset}x_iR\big)^{\un} = \big( 0 :
\sum_{i=1}^tx_iR\big)$, we have $h_{\emptyset} \in
\Gamma_{\m}(R)$, so that $h_{\emptyset}^{Q_1} = 0$.

Now suppose that $1 \le r < t$, and assume that ($\dagger$) has
been proved for $|\L| < r$. That assumption and $(*_0)$ mean that
\[
{\textstyle \sum_{\L \subsetneq [1,t], |\L| \ge r} x_{\L}^{pq-q}
h_{\L}^q \in \big(\sum_{i=1}^tx_iR\big)\bp{pq}.} \tag{$*_r$}
\]
To prove ($\dagger$) for $\L \subsetneq [1,t]$ with $|\L| = r$,
there is no loss of generality in our assuming that $\L = [1,r]$.
For every $\L' \subsetneq [1,t]$ with $|\L'| \ge r$ but $\L' \neq
[1,r]$, we have $x_{\L'}^{pq-q} h_{\L'}^q \in \sum_{i = r+1}^t
x_i^{pq-q}R$. Therefore, by $(*_r)$, we have
\[
{\textstyle x_{[1,r]}^{pq-q} h_{[1,r]}^q \in \sum_{i = 1}^t
x_i^{pq}R + \sum_{i = r+1}^t x_i^{pq-q}R = \sum_{i = 1}^r
x_i^{pq}R + \sum_{i = r+1}^t x_i^{pq-q}R.}
\]
Since $h_{[1,r]}^q \in \big(\sum_{i=1}^r x_i^qR\big)^{\un}
\subseteq \big(\big(\sum_{i=1}^r x_i^qR\big):x_t\big)$, it follows
that $${\textstyle x_{[1,r]}^{pq-q} h_{[1,r]}^q x_t \in
x_{[1,r]}^{pq-q} \big(\sum_{i = 1}^r x_i^{q}R \big) \subseteq
\sum_{i = 1}^r x_i^{pq}R}\mbox{;}$$ this implies that
\[
{\textstyle x_{[1,r]}^{pq-q} h_{[1,r]}^q \in
\big(\big(\sum_{i=1}^r x_i^{pq}R\big):x_t\big) = \big(\sum_{i=1}^r
x_i^{pq}R\big)^{\un}.}
\]
Thus $x_{[1,r]}^{pq-q} h_{[1,r]}^q \in \big(\sum_{i=1}^r
x_i^{pq}R\big)^{\un} \bigcap \big(\sum_{i = 1}^r x_i^{pq}R +
\sum_{i = r+1}^t x_i^{pq-q}R\big)$, and this is equal to
$\sum_{i=1}^r x_i^{pq}R$ by Theorem~\ref{Hu}(ii); therefore
$h_{[1,r]}^q \in \big(\sum_{i=1}^r x_i^{q}R\big)^{\lim}$.
Therefore $h_{[1,r]}^{Q_1} \in \big(\sum_{i=1}^r
x_i^{Q_1}R\big)^{\lim}$ by \ref{fixedQ}, so that
$x_{[1,r]}^{pQ_1-Q_1} h_{[1,r]}^{Q_1} \in \sum_{i=1}^r
x_i^{pQ_1}R$ by Corollary~\ref{br.21}(i). This concludes the
inductive step in the proof of ($\dagger$) and so it follows that
$y^{pQ_1} \in \big(\sum_{i=1}^t x_iR\big)\bp{pQ_1}$. This is
enough to complete the proof that $\big(\big(\sum_{i=1}^t
x_iR\big)^F\big)\bp {Q_0} = \big(\sum_{i=1}^t x_iR\big)\bp {Q_0}$
for all systems of parameters $\ul x = x_1, \dotsc, x_t$ of $R$
that are unconditioned strong $d$-sequences.

Now let $h$ be the integer of \ref{br.1b}(ii) and let $Q_4$ be a
power of $p$ with $Q_4 \ge h$. Also set $Q = Q_4 Q_0$. Let $y_1,
\ldots, y_t$ be an arbitrary system of parameters of $R$. By
Theorem \ref{br.1b}, the system of parameters $y_1^{Q_4}, \dotsc,
y_t^{Q_4}$ is an unconditioned strong $d$-sequence. Therefore, by
the first part of the proof, \begin{align*} {\textstyle
\big(\big(\sum_{i=1}^t y_iR\big)^F\big)\bp{Q_4Q_0}} & =
{\textstyle \big(\big(\big(\sum_{i=1}^t
y_iR\big)^F\big)\bp{Q_4}\big)\bp{Q_0}} \subseteq {\textstyle
\big(\big(\sum_{i=1}^t y_i^{Q_4}R\big)^F\big)\bp{Q_0}} \\ & =
{\textstyle \big(\sum_{i=1}^t y_i^{Q_4}R\big)\bp{Q_0} =
\big(\sum_{i=1}^t y_iR\big)\bp{Q_4Q_0}.} \end{align*} Thus we have
shown that ${\textstyle \big(\big(\sum_{i=1}^t x_iR\big)^F\big)\bp
{Q} = \big(\sum_{i=1}^t x_iR\big)\bp {Q}}$ for all systems of
parameters $x_1, \dotsc, x_t$ of $R$.

Finally, let $j \in \{0,\ldots, t-1\}.$ For every $n \in \mathbb
N_+$, we can apply what we have just proved to the system of
parameters $x_1, \ldots x_j,x_{j+1}^n, \ldots, x_t^n$. Thus
\begin{align*}
{\textstyle \big(\big(\sum_{i=1}^j x_iR\big)^F\big)\bp {Q}} &
\subseteq
  {\textstyle \bigcap_{n \in \mathbb N_+}
  \big(\big(\sum_{i=1}^j x_iR + \sum_{i=j+1}^t
  x_i^nR\big)^F\big)\bp{Q}}
  \\ &= {\textstyle \bigcap_{n \in \mathbb N_+}
  \big(\sum_{i=1}^j x_iR + \sum_{i=j+1}^t x_i^nR\big)\bp{Q}}\\
& = {\textstyle \bigcap_{n \in \mathbb N_+} \big(\sum_{i=1}^j
x_i^{Q}R + \sum_{i=j+1}^t x_i^{nQ}R\big)
      =\big(\sum_{i=1}^j x_iR\big)\bp {Q}}
\end{align*}
by Krull's Intersection Theorem. Now the proof is complete.
\end{proof}

\begin{corollary}
\label{mr.11} Suppose that $(R,\m)$ is a formally catenary local
ring all of whose formal fibres are Cohen--Macaulay. (These
hypotheses would be satisfied if $R$ was an excellent local ring.)
Assume further that $R$ is equidimensional, of prime
characteristic $p$ and of dimension $2$. Then there exists $Q$
such that $\big(\big(\sum_{i=1}^l x_iR\big)^F\big)\bp {Q} =
\big(\sum_{i=1}^l x_iR\big)\bp {Q}$ for all subsystems of
parameters $\ul x = x_1, \ldots,x_l$ (where $l \le 2$, of course)
of $R$.
\end{corollary}

\begin{proof} The hypotheses about $R$ are all inherited by
$R/\sqrt{0}$, and so, in view of Lemma \ref{nil}, we can assume
that $R$ is reduced. But then $R$ is Cohen--Macaulay on the
punctured spectrum, and so is a generalized Cohen--Macaulay local
ring by \ref{br.2}. The result now follows from
Theorem~\ref{dplus-Q}(ii).
\end{proof}

\end{document}